\documentclass{article}
\usepackage{amsmath, amssymb, graphics}

\newcommand{\mathsym}[1]{{}}

\begin{document}

\title{ { }  { }Non-unique factorizations, land surveying and electricity. { } { } { } { } { } { } { } { } { } { } { } { } { } { } { } { } { } {
} { } { } { } { }}
\author{Jan Sliwa}
\date{}
\maketitle

\pmb{ Abstract.} Non-unique factorizations theory, which started in algebraic number theory, over the years has expanded into several areas of mathematics.
Here, we propose yet another branching. We show that some concepts of factorizations theory, such as half factorial and weak half factorial properties
can be translated via \textit{Cayley graphs }into graph theory. It is proved, that subset \textit{S} of abelian group \textit{G} is half factorial,
if and only if the Cayley digraph \textit{Cay}(\textit{G}; \textit{S}) is geodetical, e.g., simple paths connecting a fixed pair of vertices have
the same length. Further, it is shown that the voltage digraph naturally arising from subset \textit{S} of group \textit{G} satisfies Kirchoff's
Voltage Law exactly when \textit{ S} is weakly half factorial.  In the concluding remarks, some loosely formulated ideas for further research are
presented.\\

\pmb{1. Factorization in algebraic number fields. } Let \textit{K} be an algebraic number field, \(R_K\) its ring of integers and \textit{H}(\textit{K}) the class group. An integer in \textit{ K} is said to be irreducible, if it cannot be written as a product of two non-unit elements of \(R_K\).
{ }\\
Let \textit{a} $\in $ \(R_K\) and { }\\
(1) { } { } { }\textit{ a} = \(a_1\)\(a_{2 }\)\(._{\text{...}}\) .\(a_k\) { } { } \\
be its factorization into irreducible integers. \\
If all factorizations of \textit{a} are equivalent, we say that \textit{a} has a unique factorization (is factorial) and write \textit{a} $\in
$ \textit{UF}.\\
The number of irreducible factors in (1), is called the length of factorization, and by \textit{L}(\textit{a}) we denote the set of different lengths
of factorizations of \textit{a}.  \\
If $|$\textit{L}(\textit{a})$|$ = 1 then \textit{a} is said to be half factorial (\textit{ a }$\in $ \textit{HF)}. If each element of subset \textit{ R }$\subset $ \(R_K\) is half factorial, then we call \textit{ R} to be half factorial itself and write
\textit{R }$\in $ \textit{HF}. { }\\
For \textit{ m }$\geqslant $ 1 we denote by \(G_m\)(\textit{K}) { }the set of all integers in \(R_K\) which have exactly \textit{ m} distinct lengths
of factorizations. \\
It is well known that \(R_K\) is \textit{UF} exactly when the class group of \textit{K} is trivial, and it was shown by L. Carlitz ([3]) that \(R_K\)
$\in $ \textit{HF} (or \(R_K\) = \(G_1\)(\textit{K}) ) if and only if \textit{ H}(\textit{K}) has 2 elements. Many of the factorizational properties
of { }integers in \textit{ K} could be expressed in terms of combinatorial structure of the class group \textit{H}(\textit{K}). This observation
had many implications and led to the development of a whole new area of research. \\
\\
\pmb{{ } { } { }2. Finite abelian groups.} Let \textit{G} be a finite abelian group written multiplicatively, and \textit{e} it's unit. If $\{$\(g_1\), \(g_2\),
..., \(g_k\)$\}$ is a \textit{ k}-tuple of non-unit elements of \textit{G} and \\
(2) { } { } { } { } { }\(\underset{i=1}{\overset{k}{\Pi }}\)\(g_i\) = \textit{ e}\\
then $\{$\(g_1\), \(g_2\), ... , \(g_k\)$\}$ will be called a block, and the block will be called minimal (or atom), if from { }\(\underset{i=1}{\overset{k}{\Pi
}}\)\(g_i^{m_i}\) = \textit{ e} where { }\(m _i\) $\in $ $\{$0,1$\}$ for 1$\leqslant $ \textit{i} $\leqslant $\textit{ k} { }it follows { }that
all \(m_i\)'s are ones or all of them are zeros. { }\\
We will denote by $\mathcal{A}$ = \(\mathcal{A}_{G }\) the set of all atoms in \textit{G} and if we restrict { }\(g_i\)'s to a subset \textit{S}
of \textit{ G}, then it will be denoted by \(\mathcal{A}_{G, S}\). \\
For a block  $\beta $ { }= $\{$\(g_1\),\( g_2\), ..., \(g_k\)$\}$ $\in $ $\mathcal{A}$ { } the set of distinct \(g_i\)'s will be called support of
$\beta $ and denoted by supp($\beta $). \\
For $\beta $ $\in $ $\mathcal{A}$ we will define \textit{ c}($\beta $), the cross number of $\beta $ to be \\  \\
(3){ } { } { } { } { } \(\text{      }\sum _{\text{\textit{$i$}}=1}^{\text{\textit{$k$}}} \frac{1}{\text{\textit{ord}}\left(\text{\textit{$g$}}_{\text{\textit{$i$}}}\right)}\)\\  \\
where \textit{ord}(\textit{g}) denotes the order of \textit{g }in G. \\
If $\beta $ is an atom satisfying { }\textit{c}($\beta $) = 1 { }then it will be called \textit{C}-\textit{atom}.\\
If \textit{S} is a subset of \textit{G,} and if each atom with support in \textit{S} is a \textit{C}-atom, then we will say that \textit{S}
has { }(\textit{C}) { }property (or satisfies condition (\textit{C})) and in this case we will write \textit{S }$\in $\textit{ C}. \\
We will denote by $\mu $(\textit{G}) the maximal cardinality of \textit{S} $\subset $ \textit{G} satisfying condition (\textit{C}), and by \textit{t}(\textit{G}) the minimal number of summands in \textit{ G }=\( \underset{i=1}{\overset{t}{\cup }}\)\(S_i\), where all { }\(S_i\)'s have property
(\textit{C}). \\
\\
\pmb{{ } { } { }3. Half factorial and weakly half factorial properties.} { }For subset \textit{ S} of the class group \textit{ H}(\textit{K}) we denote by
\(R_{K,S}\) the set of integers \textit{a} $\in $ \(R_K\), such that all prime ideals dividing the principal ideal (\textit{a}) = \textit{a}\(R_K\)
correspond to elements of \textit{ S} under the natural mapping of \(\text{ideals}\) into \textit{H}(\textit{K}). \\
The connection between the factorizations of elements in \(R_{K,S}\) and (\textit{C}) property of { }\textit{S} was first observed in [41]. \\
It was noticed there that \(R_{K, S}\) is \textit{ HF} (or \(R_{K, S}\) $\subset $ \(G_1\)) if and only if \textit{ S} satisfies condition (\textit{C}). \\
The relationship between the cross number of atomic blocs and factorizations in semigroups and monoids was later observed by L. Skula ([40]) and
A. Zaks ([45]). \\
In recent terminology, the subset of the class group having property (\textit{C}) is also said to be half factorial.\\
In ([43]) the structure of \(G_m\)(\textit{K}) was determined and it was shown that if the class number of K is greater than 2, then the set $\{$\textit{a} $\in $ \(G_m\)(K): { }\(\text{Norm}_{K/Q}\)(\textit{ a}) $\leqslant $ \textit{x}$\}$ (if non-empty) has an asymptotic density of the type:\\
\\
 { } { } { } { } { } \(G_m\)(\textit{x}) = (\textit{C} +\textit{ O}(1))\(x(\log  \log  x)^A\) /(\(\text{\textit{$\log  x$}})^B\), \\
\\
where A and C are some combinatorial constants that depend on \textit{m} and \textit{ H}(\textit{ K}) and { }\(B =1-\frac{\mu (H(K))}{|H(K)|}\).
\\
Similar asymptotic evaluations were shown to be true for set of positive rational integers $\leqslant $ \textit{x }which have \textit{m} lengths
of factorization in \(R_K\). \\
In the subsequent paper ([44]) it was proved that if { }$|$\textit{H}(\textit{K})$|$ $\geqslant $ 3 then in the union\\
\\
 { } { } { } { } { } \(R_{K }\)= \(\underset{m=1}{\overset{\infty }{\cup }}\)\(G_m\)(\textit{K}), \\
\hspace*{3.5ex} \\
all summands are non-empty and a weaker than (\textit{C}) property was introduced. Namely, we said there that subset \textit{ S} of abelian group
\textit{ G} has property (\(C_0\)) (or is weakly half factorial - using recent terminology) if for each atom $\beta $ with support in \textit{ S,}
the\textit{  }cross number \textit{ c}($\beta $) is an integer. That leads to the definitions of constants \(\mu _0\)(\textit{G}) and \(t_0\)(\textit{G}). \\As { }\(\mu _0\)(\textit{G}) $\geqslant $ \(\mu \)(\textit{G}) and \(\mu _0\)(\textit{G}) is easier to evaluate, { }investigating the
sets with property (\(C_0\)) gives some nontrivial upper bound for $\mu $. { } { } { } { }\\
The results on structure of \(G_m\)(\textit{K}) that were proven in [43], were stated for the ring of algebraic integers, but it is easily seen
that they can be formulated and proven in the more general setting of block semigroup over \textit{H}(\textit{K}). \\
Actually, papers [41]-[44] already contained the algebraic, analytical and combinatorial elements of a later developed theory. \\
\\
\pmb{{ } { } { }4. Non-unique factorizations theory. }The result of L. Carlitz ([3]), a series of papers by W. Narkiewicz on densities of integers with unique/non-unique
factorizations { }in quadratic fields ([23]-[28]), { }and { }the above mentioned papers by the current author ([41]-[44]), gave birth to the non-unique
factorizations theory in algebraic number fields. \\This, with the subsequent result of L. Skula ([40]) and A. Zaks ([45]) on factorizations in semigroups
and monoids, attracted a lot of attention a few years later. A vast volume of papers was published. Most of them deal with factorizations in Dedekind
domains, Krull monoids, asymptotic density of elements possessing specific factorizational properties, the structure of the set lengths of factorizations,
and the values of some combinatorial constants defined for abelian groups that are related to factorizations in block semigroups. { }\\
\\
It is not my intention here, and it would be difficult, to give an overview of the current state of the theory of non-unique factorizations in a
short note. Almost randomly and arbitrarily, I will refer the reader to recent books, proceedings and papers by the most active authors in the area
: { }D. D. Anderson ([1]), S. T. Chapman ([4]-[6]), W. Gao ([8],[9]), A. Geroldinger ([10]-[14]), F. Halter-Koch ([14],[16],[17]), W. Hassler ([18]),
J. Kaczorowski ([14],[19]-[21]), U. Krause ([22]), W. Narkiewicz ([29]-[32]), A. Plagne ([33],[34]), M. Radziejewski ([35]-[38]) and W. A. Schmid
([37],[38],[39]).\\
\\
As it can be seen in the literature quoted above, non-unique factorizations theory has been applied to various areas of mathematics using a broad
range of techniques. In the current note, we propose branching some of the factorizational concepts into graph theory. { }In recent decades, graph
theory has advanced progress in many areas of mathematics, and was itself helped by such intrusions. \\ \\
\pmb{ { } { } { }5. Basic terms in digraphs and graphs. }We adopt the standard graph theory notation and terminology. For definition of terms not defined here, and
quoted results, the reader is referred to [2], [7] or [15]. \\
Let \textit{D} = (\textit{V}, \textit{A}) be a digraph (directed graph) with finite set of vertices \textit{V }= \textit{V}(\textit{D}) { }and arcs \textit{ A }= \textit{A}(\textit{D}). We consider here only the digraphs with no loops and with no multiple arcs. That means we
can consider each arc \textit{ a} $\in $ \textit{A}(\textit{D}) an ordered pair of distinct vertices, and write { }\textit{ a = }(\textit{x},
\textit{y}), where \textit{ x}, \textit{y} $\in $ \textit{V. }We will call \textit{x }and\textit{ y the beginning} and \textit{ending} of arc
\textit{a} and write \textit{ x }= \textit{init}(\textit{a}), \textit{y }= \textit{fin}(\textit{a}). A \textit{walk} \textit{W} in \textit{D} is an alternating sequence \\ 
  \textit{  W }: \(x_1\)\(a_1\)\(x_2\)\(a_2\)\(x_3\) . . . \(x_{k-1}\)\(a_{k-1}\)\(x_k\)\(a_k\)\(x_{k+1}\) \\
of vertices \(x_i\) and arcs \(a_i\) such that the beginning of arc \(a_i\) is vertex { }\(x_i\) and the ending of \(a_i\) is \(x_{i+1}\). As there
are no multiple arcs, we can just write \\ \textit{ W }: \(a_1\)\(a_2\) . . . \(a_k\)\(.\)  \\A walk is called a \textit{path} if all of its vertices
are distinct, with the possible exception of the first and last. In this last case, e.g., if  \(x_1\) = \(x_{k+1}\) the path will be called a
\textit{cycle}. \\
If for path \textit{ P} there is { }\(x_1\) = \textit{x} and \(x_{k+1}\) = \textit{y,} then \textit{P} will be called (\textit{x}, \textit{y})-\textit{path}. In particular, for any vertex \textit{x} belonging to a cycle, we can consider the cycle to be (\textit{x}, \textit{x})-path.
As for the arcs, we can talk of beginning and ending vertices of the walk, and define functions \textit{fin} and \textit{init} from a walk set
into vertices. \\
In case the digraph \textit{D} is \textit{symmetrical}, e.g., if (\textit{x}, \textit{y}) $\in $ \textit{D} always implies (\textit{y}, \textit{ x}) $\in $ \textit{ D}, we can ignore the direction of the arcs, call them edges and consider \textit{ D} to be a \textit{ graph}. So
graphs, in a sense, could be considered a special case of digraphs. Even if the digraph is not symmetrical, we can, by ignoring the direction of
the arcs, talk about the \textit{underlying graph} of digraph D. \\
All the terms that we will introduce for the digraphs, can be applied to graphs. \\
If \textit{D} = (\textit{V}, \textit{A}) is a digraph, then any function $\psi $ :\textit{A} $\longrightarrow $ \(R^+\) we will call \textit{weight}. Its reciprocal, { }$\lambda $\textit{  }= 1/$\psi $ { }will be called the \textit{length,} and the pair (\textit{D}, $\psi $) forms the
\textit{weighted digraph}. We will require that if { }(\textit{x}, \textit{y}) and (\textit{y}, \textit{x}) are arcs of the digraph, then they
both have the same weight, so that there will be no confusion when going from the digraph to its underlying graph. { } \\Weight $\psi $, and length
$\lambda $ could be additively extended to any walk \textit{ W}: \(a_1\)\(a_2\) ... \(a_k\) by setting $\psi $(\textit{W}) = \(\sum _{i=1}^k\)$\psi
$(\(a_i\)) { }and { }$\lambda $(\textit{W}) = \(\sum _{i=1}^k\)$\lambda $(\(a_i\)). Usually, by the length of the walk, the number of arcs connecting
its beginning and ending is understood. In our terminology this corresponds to the particular choice of weight, { }namely $\psi $ $\equiv $ 1.\\
\\
If C is a finite set of colors, then any function \textit{ c}:\textit{A }$\longrightarrow $ \textit{C} is called \textit{ arc coloring} and the
digraph is said to be \textit{ k-colored} if $|$\textit{C}$|$ = \textit{ k}. If none of the adjacent arcs are colored with the same color, then
coloring is said to be \textit{ proper}. \\
The set of all paths (including closed paths e.g., cycles) in \textit{ D} will be denoted by { }$\mathcal{P}$ = \( \mathcal{P}_D\) and the set of
all (\textit{init}(\textit{P}), \textit{fin}(\textit{P}))-paths { }we will denote by \(\mathcal{P}_D\)(\textit{P}). The elements of\textit{
 } \(\mathcal{P}_D\)(\textit{P}) { }will be referred to as \textit{ P}-paths. \\
\\
\\
\pmb{{ } { } { }6. Some more definitions in graph and digraph theory.} \\
\pmb{{ } { } { }\(6.1 \text{ Cayley} \text{ digraphs}.\)} Let \textit{ G} be a group and \textit{ S} its subset. The \textit{ Cayley} \textit{ digraph,} \textit{
Cay}(\textit{G}; \textit{S}) is a digraph whose vertices are identified with elements of G\textit{ , }and the arcs set consists of ordered pairs
{ }(\textit{g}, \textit{gs}), where \textit{ g} $\in $ \textit{G} and \textit{ s }$\in $ \textit{S}. We do not require that \textit{S} generate
\textit{G}, hence the resulting digraph is not necessarily connected\textit{ .} In case \textit{ S} is symmetrical, e.g., { }\textit{S} = \(S^{-1}\)
we can also talk about the \textit{Cayley graph}. Cayley digraphs are also called \textit{Cayley color digraphs,} as each of its arcs  (\textit{g}, \textit{gs}) could be considered colored by an element \textit{ s} of \textit{ S}.\\
\hspace*{0.5ex} \\
\pmb{{ } { } { }6.2 Geodetical Digraphs.} Let \textit{ D} = (\textit{V}, \textit{A}) be a $\psi $-weighted digraph. Path \textit{ P} is called \textit{geodesic}
if it is the shortest path connecting \textit{ init}(\textit{P}) and \textit{ fin}(\textit{P}). If each path connecting \textit{ init}(\textit{P}) and \textit{ fin}(\textit{P}) is geodesic, then \textit{ P} is called \textit{ geodetical}. In other words, \textit{ P} is geodetical if all
(\textit{init}(\textit{P}), \textit{fin}(\textit{P}))-paths have the same length or { }$|\{\lambda $(\textit{Q}): \textit{ Q }$\in $ \(\mathcal{P}_D\)(\textit{P})$\}|$ = 1. \\
If each path in \textit{ D} is geodetical, then the digraph itself will be called \textit{geodetical. }(We use the term \textit{geodetical}, as both
\textit{geodesical} and \textit{gedotetic} are already taken for describing graphs with different properties). \\
\hspace*{0.5ex} \\
\pmb{{ } { } { }6.3  Voltage digraphs.} If \textit{ D} = (\textit{V}, \textit{A}) { }is a $\psi $-weighted digraph, \textit{ H} a group and $\phi $:\textit{ A} $\longrightarrow $ \textit{H} is mapping from arcs set into \textit{ H,} then the triple { }(\textit{D}, $\phi $, \textit{H}) { }is called
a \textit{voltage digraph}. Voltage digraphs and graphs (defined similarly, with additional requirements that { }$\phi $(\(a^{-1}\)) = \( \phi (a)^{-1}\)
) are powerful tools for constructing large graphs, called lifts with prescribed properties as covering spaces of small base graphs. { }\\
If for a closed path \textit{ P}: \(a_1\)\(a_2\) ... \(a_k\) { }there is { }\(\prod _{i=1}^k \)\(\phi \left(a_i\right.\)) = \textit{ e, }where \textit{
e } is the unit element of \textit{ H}, then we say that path \textit{ P} satisfies \textit{ Kirchoff's Voltage Law} (in short \textit{ P} $\in $
\textit{ KVL)}. If each closed path satisfies \textit{ KVL,} then we say that the voltage digraph { }(\textit{D}, $\phi $, \textit{H}) { }satisfies
Kirchoff's Voltage Law. \\
\hspace*{1.5ex} \\
\pmb{{ } { } { }7. Cayley digraphs based on half factorial and weakly half factorial sets. }\\Let \textit{G} be an abelian group, \textit{ S} its subset and
\textit{D} = \textit{Cay}(\textit{G}; \textit{S}) be the Cayley digraph. For arc \textit{ a} = (\textit{g}, \textit{gs}) with \textit{ g }$\in
$ \textit{G} and \textit{ s }$\in $ \textit{S} we will set $\psi $(\textit{a}) = \textit{ord}(\textit{s}) (and hence $\lambda $(\textit{a})=\(\frac{1}{\text{\textit{ord}}(s)}\)).\\
\\
\pmb{Theorem 1. }Subset \textit{S} of \textit{G} is half factorial if and only if digraph \textit{Cay}(\textit{G}; \textit{S}) is geodetical.
{ }\\
\\
Proof. Let \textit{D} = \textit{Cay}(\textit{G}; \textit{S}), \textit{ x, y} $\in $ \textit{V}(\textit{D}) and let \(P_1\) : \(a_1^{(1)}\)\(a_2^{(1)}\)
... \(a_m^{(1) }\) and \(P_2\) : \(a_1^{(2)}\)\(a_2^{(2)}\) ... \(a_n^{(2)}\) be { }(\textit{x,y})-paths, e.g. beginnings of arcs \(a_1^{(1)}\),
\(a_1^{(2)}\) and endings of \(a_m^{(1)}\), \(a_n^{(2)}\) coincide. { }\textit{  { }\\
}Let { }\(a_i^{(1)}\) = (\(g_i^{(1)}\), \(g_i^{(1)}\)\(s_i^{(1)}\)), \(a_j^{(2)}\) = (\(g_j^{(2)}\), \(g_j^{(2)}\)\(s_j^{(2)}\)) where \(g_i^{(1)}\),
\(g_j^{(2)}\) $\in $ \textit{ G} and { }\(s_i^{(1)}\), \(s_j^{(2)}\) $\in $ \textit{ S} for \textit{ 1}$\leqslant $\textit{  i $\leqslant $ m, 1$\leqslant
$ j $\leqslant $ n}. Without loss of generality we can assume that \textit{ x = e, }unit element { }of \textit{ G}.\\
The fact that \(P_1\) and \(P_2\) are (\textit{x, y})-paths could be restated as \\
 { } { } { } { } { } { } { }\(s_1^{(1)}\)\(s_2^{(1)}\) ... \(s_m^{(1) }\) =  \(s_1^{(2)}\)\(s_2^{(2)}\) ... \(s_n^{(2) }\).\\
Choose \(s_1^{(3)}\), \(s_2^{(3)}\), ... , \(s_k^{(3)}\) { }such that\\
\\
\hspace*{1.5ex} (4) { } { } { }\(\prod _{i=1}^m \)\(s_i^{(1)}\) \(\prod _{t=1}^k \)\(s_t^{(3)}\) = { }\(\prod _{i=1}^n \)\(s_i^{(2)}\) \(\prod _{t=1}^k
\)\(s_t^{(3)}\) = \textit{ e { }\\
\hspace*{1.5ex} \\
}and that \\
\\
 { } { } { } { } { } { }\(\prod _{i=1}^m \)\(s_i^{(1)}\) \(\prod _{t\in \{1,2, \text{...}, k\}\left\backslash t_{0\text{   }}\right.} \)\(s_t^{(3)}\)
\\
\hspace*{1.5ex} \\
is not equal to \textit{ e} for any \(t_0\)$\in $ $\{$1,2, ..., k$\}$. Let us write the left and right sides of (4) as the products of atomic blocks
\(\beta _1^{(1)}\)\(\beta _2^{(1)}\) ... \(\beta _{l_1}^{(1)}\) and \(\beta _1^{(2)}\)\(\beta _2^{(2)}\) ... \(\beta _{l_2}^{(2)}\), respectively.
From \textit{ S }$\in $ \textit{ HF} it follows that \(l_1\) =\(\text{  }l_2\). Hence \\
\hspace*{1.ex} \\
\hspace*{6.ex} $\lambda $(\(P_1\)) = \(l_1\) - \(\sum _{t=1}^k \)\(\frac{1}{\text{\textit{ord}}\left(s_t^{(3)}\right)}\) = \(l_2\) - \(\sum _{t=1}^k
\)\(\frac{1}{\text{\textit{ord}}\left(s_t^{(3)}\right)}\) = $\lambda $(\(P_2\)) . \\
Conversely, let us assume that \textit{Cay}(\textit{G}; \textit{S}) is geodetical with weight defined by { }$\psi $((\textit{g}, \textit{gs}))
= \textit{ ord}(\textit{s}). Let { }$\beta $ : \(s_1\)\(s_2\) ... \(s_k\) { }be an arbitrary atom with support in \textit{ S} and \(P_1\) be the
(\textit{e,e})-path corresponding to $\beta $, e.g., { }\(P_1\) : (\textit{e}, \(s_1\))(\(s_1\), \(s_1\)\(s_2\)) ... (\(s_1\)\(s_2\) ... \(s_{k-1}\),
\(s_1\)\(s_2\) ... \(s_k\)). For \(P_2\) we choose (\textit{e, e})-path corresponding to the atomic block \( s_1\)\(s_1\) ... \(s_1\) { }(\(s_1\)
taken \textit{ord}(\(s_1\)) times). Obviously\textit{  }$\lambda $(\(P_2\)) = 1, hence $\lambda $(\(P_1\)) = 1. This means that { }each atomic block
with support in \textit{ S} is a \textit{ C}-atom. That concludes proof of Theorem 1.\\
\hspace*{1.5ex} \\
Let \textit{D} be a Cayley digraph \textit{Cay}(\textit{G}; \textit{S}) and \(C_N\) the group of \textit{N}-th root of unity, where \textit{N} is the exponent of \textit{G. }We will define $\phi $ by setting $\phi $((\textit{g}, \textit{gs})) = \textit{ exp}(\(\frac{2\text{$\pi $i}}{\text{\textit{ord}}(s)}\)).
Triple (\textit{G}, $\phi $, \(C_N\)) constitutes a voltage digraph. \\
\\
\pmb{ Theorem 2}. Subset \textit{S} of group \textit{G} is weakly half factorial if and only if the voltage digraph (\textit{Cay}(\textit{G}; \textit{S}), $\phi $, \(C_N\)) { }satisfies Kirchoff's Voltage Law.\\
\\
Proof. Let \textit{ P }: \(a_1\)\(a_2\) ... \(a_k\) be a closed path where \(a_i\) = (\(g_i\), \(g_i\)\(s_i\)) and \(g_{i }\)$\in $ \textit{ G},
\(s_i\)$\in $\textit{ S} for \textit{ i = 1, 2, ... , k}. We can assume that the starting and ending vertices of P coincide with \textit{e}, the
unit element of \textit{ G}. The fact that \textit{ P} is a closed path translates to equality\\
\hspace*{1.ex} (5) { } { } { } { } { }\(\underset{i=1}{\overset{k}{ \prod }}\) \(s_i\) = \textit{ e.\\
}If \textit{ S} $\in $ \textit{WHF} { }then { } \(\sum _{i=1}^k \)\(\frac{1}{\text{\textit{or}}d\left(s_i\right)}\) is an integer and this implies \\
{ }$\phi $(\textit{P}) = \(\prod _{i=1}^k \)$\phi $(\(a_i\)) = \textit{ exp}(2$\pi $\textit{i}\(\sum _{i=1}^k \)\(\frac{1}{\text{\textit{ord}}\left(s_i\right)}\))
= 1.\\
Now assume that the voltage digraph (\textit{Cay}(\textit{G}; \textit{S}), $\phi $, \(C_N\)) where \(C_N\) and $\phi $ are defined as above, satisfies
Kirchoff's Voltage Law, and let $\beta $: \(s_1\)\(s_2\) ... \(s_k\) be a block with support in \textit{S}. Path \textit{P} corresponding to $\beta
$ is a closed (\textit{e}, \textit{e})-path, so KVL implies that \textit{ exp}(2$\pi $\textit{i}\(\sum _{i=1}^k \)\(\frac{1}{\text{\textit{ord}}\left(s_i\right)}\))
= 1. Hence $\beta $ $\in $ (\(C_0\)). This proves that \textit{ S }$\in $ \textit{WHF}.\\
\\
\pmb{{ } { } { }8. Concluding remarks. }We will sometimes call geodetical digraphs to be half factorial and voltage digraphs satisfying KVL - weakly half factorial.
Both of these terms could be extended to graphs, hence in what follows, we will mean both graphs and digraphs, even if we specifically mention one
or the other. { }\\The following remarks suggest some further possibilities for applying factorizational concepts in the field of graph theory. { }
\\
\\
\pmb{{ } { } { }8.1 Nonabelian groups. }The half factorial concept has its origin in factorizations of algebraic integers and through property (\textit{C}) it is
being transferred to subsets of the class group, and from there to any abelian group. Theorem 1 and Theorem 2 translate this concept and weakly half
factoriality to digraphs. Both of them could be used in reverse, as definitions of half factorial and weakly half factorial properties of the subset
of the group. Notice, that the requirement that the group is abelian is dropped. This gives reason to study half and weakly half factoriality for
nonabelian groups. \\
\\
\pmb{{ } { } { }8.2 Number of lengths of paths.} Let $\mathcal{P}$ be the set of paths in the graph/digraph and let \textit{P} $\in $ $\mathcal{P}$. By
{ }$\mathcal{L}$(\textit{P}) { }we denote the set of distinct lengths of all \textit{P}-paths, e.g., { }$\mathcal{L}$(\textit{P}) = $\{\lambda
$(\textit{Q}): \textit{ Q} $\in $ $\mathcal{P}$(P)$\}$. \\
Hence, it could be said that \textit{P} is geodetical if and only if  $\#\{\mathcal{L}$(\textit{P})$\}$ = 1,\textit{  }and the graph/digraph
is geodetical if this is true for each element of { }$\mathcal{P}$. This leads us to the definition of \textit{ m}-geodetical paths and sets of paths.\\
Namely, we will say that \textit{P} $\in $ $\mathcal{P}$ is \textit{ m-geodetical} if { } $\#\{\mathcal{L}$(\textit{P})$\}$ = \textit{ m,} and
by \(G_m\) we will denote $\{$ \textit{P} $\in $ $\mathcal{P}$: { }$\#\{\mathcal{L}$(\textit{P})$\}$ = \textit{ m }$\}$. \\
It would be interesting
to know under what conditions { }\(G_m\)(\textit{D}) $\neq $ $\{\varnothing \}$ implies that { }\(G_i\)(\textit{D}) $\neq $ $\{\varnothing \}$
for \textit{ i }$\leqslant $ \textit{ m, }and what is the largest value of \textit{m}  such that \(G_m\)(\textit{D}) $\neq $ $\{\varnothing \}$.\textit{ \\
}\\
The uniqueness of the factorizations concept, might be in a natural way transferred to the graph theory by introducing the unique path property.
Namely, we will say that \textit{ P} $\in $ $\mathcal{P}$ has \textit{unique path property } (\textit{P} $\in $ \textit{UP}) { }if \textit{P} is the unique element in { }$\mathcal{P}$(\textit{P}). { }It is well known that only connected \textit{UP} graphs are trees, and so generally
it could be said that the graph is \textit{ UP} if and only if it is a forest. The characterization of \textit{UP} digraphs is slightly more complex. \\
Obviously \textit{UP} digraphs/graphs are geodetical. \\In the case of factorizations of algebraic integers, the structure of the class group provides
the answer to the question how far the ring of integers is from being \textit{UF}. Here, we can suspect that the complexity of $\{\mathcal{L}$(P):
P $\in $ $\mathcal{P}\}$ depends on how different is the graph from being a tree or acyclic. The next observation gives some weight to this speculative
statement. { }\\
\\
Let \(R^V\) and \(R^A\) be the finite dimensional vector spaces of functions \textit{p}:\textit{V }$\longrightarrow $ \textit{R} and
\textit{g}:\textit{A }$\longrightarrow $ \textit{R, }respectively. Let further $\delta $: \(R^V\)$\longrightarrow $ \(R^A\) be the
linear transformation defined for \textit{ p }$\in $ \(R^V\) by $\delta $\textit{p}(\textit{a}) = \textit{ p}(\textit{fin}(\textit{a})) -
\textit{p}(\textit{init}(\textit{a})). \\Subspace { }\(\mathcal{B}\) = $\delta $(\(R^V\)) { }of \( R^{A }\) is called \textit{ bond space} and
is related to \textit{cycle space}. For each element \textit{ g} $\in $ \(\mathcal{B}\) there exists \textit{ p }$\in $ \(R^V\) such that \textit{
g }= $\delta $\textit{p}, hence for each (\textit{x}, \textit{y})-path \textit{P}: \(a_1\)\(a_2\) ... \(a_k\) we have \\
\\
 { } { } { } { } { }\textit{ g}(\textit{P}) = \(\sum _{i=1}^k \)\textit{ g}(\(a_i\)) = \(\sum _{i=1}^k \)(\textit{p}(\(\text{\textit{fin}}\left(a_i\right.\))
- \textit{p}(\textit{init}(\(a_i\))) = \textit{p}(\textit{y}) - \textit{p}(\textit{x}). \\
\hspace*{5.ex} \\
Therefore, if we set \(A_g\)= $\{$\textit{a }$\in $ \textit{ A}: \textit{ g}(\textit{a}) $>$ 0$\}$, { }\(\psi _g\)(\textit{a})= 1/\textit{g}(\textit{a}), the weighted digraph \( D_g\)\textit{  }= { }(\textit{V}, \(A_g\), \(\psi _g\)) { }will be geodetical. \\
As $\mathcal{B}$ is an orthogonal complement
(with suitably selected scalar product) to cycle space, we can speculate that lesser cyclicity of the diagraph implies richer geodeticality. { }\\
\\
\pmb{{ } { } { }8.3 Colored digraphs/graphs and half factoriality.} Let \textit{G }= (\textit{V}, \textit{E) { }}be a finite $\psi $-weighted, properly
k-colored graph. By $\mu $(\textit{G}) we denote the maximal edge size of the geodetical subgraph (though not necessarily an induced subgraph)
of \textit{G} and by \textit{t}(\textit{G}) the minimal number of such subgraphs that  cover \textit{G. \\
}There is { }\\
\hspace*{0.5ex} \\
\hspace*{0.5ex} (6) { } { } { }t(G)$\leqslant $ \(\chi \)'(G), { }$\mu $(G)$\leqslant |$\textit{ E}$|$/\(\chi \)'(G)\\
\hspace*{0.5ex} \\
where \(\chi \)'(G) is the edge-chromatic number of G . \\
To show this, observe that if \textit{ C }= $\{$\(c_1\), \(c_{2 }\), ..., \(c_k\)$\}$ is the set of colors with which \textit{ G} is colored and
{ }\(E_i\) = $\{$\textit{e }$\in $ \textit{ E}: \textit{ c}(\textit{e}) =\( c_i\)$\}$ for { }\textit{ i }= 1, ... , \textit{ k} { }then { }\(E_1\),
\(E_2\), ..., { }\(E_k\) cover \textit{ E. }As the coloring is proper, each of \(E_i\) is \textit{ UP, }hence it is geodetical. Therefore { }t(\textit{G}) $\leqslant $ \textit{ k} and $\mu $(\textit{G}) $\leqslant $ { }\(\min \) $\{|$\(E_i\)$|$: 1$\leqslant $i$\leqslant $\textit{ k}$\}$. Now
by definition, { }\(\chi \)'(\textit{G}) is minimal \textit{ k}, such that \textit{ G} is properly \textit{ k}-edge colorable. As \( \chi \)'(\textit{G}) { }is equal { }\textit{$\Delta $}(\textit{G}) { }or { }\textit{$\Delta $}(\textit{G})+1, where \textit{ $\Delta $}(\textit{G}) { }is
the degree of graph \textit{G}, (6) could be restated in the form:\\
\hspace*{0.5ex} \\
 { } { } { } { } { } { }\textit{ t}(\textit{G})$\leqslant $ \textit{ $\Delta $}(\textit{G})+1 { }and $\mu $(\textit{G})$\leqslant |$\textit{ E}$|$/\textit{
$\Delta $}(\textit{G}).\\
\hspace*{1.5ex} \\
A similar evaluation could be shown to hold for digraphs.\\
\\
The observation, that the monochromatic set of edges in properly colored graphs constitutes a half factorial subgraph (with vertices set V), suggests
an alternative definition of constant $\mu $(G), more in line with the one implied by the Cayley digraphs. \\Namely, if { }\textit{S} { }is a subset
of { }\textit{C { }}and { }\(A_{S }\) = $\{$\textit{a }$\in $ \textit{A}: \textit{c}(\textit{a}) $\in $ \textit{S}$\}$, then by \(\mu ^{\star
}\)(\textit{G}) we will denote the maximal cardinality of S, such that \(A_S\) is half factorial, and by \(t^{\star }\)(\textit{G}) minimal \textit{
t}, such that \textit{ S }=\( \underset{i=1}{\overset{t}{\cup }}\)\(S_i\), where all \(S_i\)'s are geodetical.\\
Similarly, we can define the constants \( \mu _0\)(\textit{G}), \(t_0\)(\textit{G}), \(\mu _0^{\star }\)(\textit{G}) and \(t_0^{\star }\)(\textit{G}) for graphs (or digraphs). \\
\\
Hopefully, using the vast machinery of graph theory some nontrivial results in the suggested directions could be obtained and the ideas that proved
themselves fruitful in non-unique factorizations theory will repeat their success in graph theory. \\ \\
\textit{ \\
}\pmb{ References.\\
\\ 
}\\1. Anderson, D.D. (ed.), Factorization in integral domains, Lecture Notes in Pure and Appl. Math. 189, Marcel Dekker Inc., New York 1996.\\
2. Bang-Jensen, J. and Gutin, G,, Digraphs: Theory, Algorithms and Applications, Springer-Verlag, London, 2000.\\
3. Carlitz, L., A characterization of algebraic number fields with class number two, Proc. Amer. Math. Soc. 11 (1960), 391-392.\\
4. Chapman, S.T. and Smith, W.W., Factorization in Dedekind domains with finite class group, Israel J. Math. 71 (1990), 65-95.\\
5. Chapman, S.T. and Glaz, S. (eds.): { }Non-Noetherian commutative ring theory, Math. Appl. 520, 97-115, Kluwer Acad. Publ., Dordrecht 2000.\\
6. Chapman, S.T. (ed.), Arithmetical properties of commutative rings and monoids, { }Lecture Notes in Pure and Appl. Math. 241, CRC Press, 2005.\\
7. Chartrand, G., Lesniak, L., Graphs $\&$ Digraphs, Chapman $\&$ Hall, Boca Raton, 2005. { }\\
8. Gao, W. and Geroldinger, A., Half-factorial domains and half-factorial subsets of abelian groups, Houston J. Math. 24 (1998), 593-611.\\
9. Gao, W. and Geroldinger, A., Systems of sets of lengths. II, Abh. Math. Sem. Univ. Hamburg 70 (2000), 31-49.\\
10. Geroldinger, A., Uber nicht-eindeutige Zerlegungen in irreduzible Elemente, Math. Z. 197 (1988), 505-529. \\
11. Geroldinger, A., Halter-Koch, F., Non-unique Factorizations in Block Semigroups and Arithm. Appl., Math. Slovaca, 42 (1992) 641-661. \\
12. Geroldinger, A., Kaczorowski, J., Analytic and arithmetic theory of semigroups with divisor theory, J. Theorie d. Nombres Bordeaux 4 (1992),
199-238.\\
13. Geroldinger, A., Schneider, R., The cross number of finite abelian groups. II, European J. Combin. 15 (1994), 399-405.\\
14. Geroldinger, A., Halter-Koch, F., Kaczorowski, J., Non-unique factorizations in orders of global fields, J. Reine Angew. Math. 459 (1995), 89-118.\\
15. Gross, J. L., Yellen, J., Handbook of Graph Theory, Boca Raton London New York Washington, D.C., 2004 CRC Press LLC.\\
16. Halter-Koch, F., Chebotarev formations and quantitative aspects of nonunique factorizations, Acta Arith. 62 (1992), 173-206.\\
17. Halter-Koch, F. and Muller, W., Quantitative aspects of nonunique factorization: a general theory with applications to algebraic function fields,
J. Reine { } { } { } { }Angew. Math. 421 (1991), 159-188.\\
18. Hassler, W., A note on half-factorial subsets of fnite cyclic groups, Far East J. Math. Sci. (FJMS) 10 (2003), 187-197.\\
19. Kaczorowski, J., Some remarks on factorization in algebraic number fields, Acta Arith. 43 (1983), 53-68.\\
20. Kaczorowski, J., Completely irreducible numbers in algebraic number fields, Funct. Approximatio, Comment. Math. 11 (1981), 95 - 104.\\
21. Kaczorowski, J., A pure arithmetical characterization for certain fields with a given class group, Colloq. Math. 45 (1981), 327 - 330.\\
22. Krause, U., and Zahlten, C., Arithmetic in Krull monoids and the cross number of divisor class groups, Mitt. Math. Ges. Hamburg 12 (1991), 681-696.\\
23. Narkiewicz, W., On algebraic number fields with non-unique factorization. Colloq. Math. 12 1964 59--68.\\
24. Narkiewicz, W., On algebraic number fields with non-unique factorization. II. Colloq. Math. 15 1966 49--58.\\
25. Narkiewicz, W., On natural numbers having unique factorization in a quadratic number field. Acta Arith. 12 1966 1--22.\\
26. Narkiewicz, W., Factorization of natural numbers in some quadratic number fields. Colloq. Math. 16 1967 257--268\\
27. Narkiewicz, W., On natural numbers having unique factorization a quadratic number field. II. Acta Arith 13 1967/1968 123--129\\
28. Narkiewicz, W., A note on factorizations in quadratic fields. Acta Arith. 15 (1968) 19--22.\\
29. Narkiewicz, W., Finite abelian groups and factorization problems, Colloq. Math. 42 (1979), 319-330.\\
30. Narkiewicz, W., Sliwa, J., Finite abelian groups and factorization problems II, Colloq. Math. 46 (1982), 115 - 122.\\
31. Narkiewicz, W.,A note on elasticity of factorizations. J. Number Theory 51 (1995), no. 1, 46--47.\\
32. Narkiewicz, W., Elementary and analytic theory of algebraic numbers, third edition, Springer-Verlag, Berlin 2004.\\
33. Plagne, A. and Schmid, W.A., On large half-factorial sets in elementary p-groups: Maximal cardinality and structural characterization, Isreal
J. Math. 145 (2005), 285-310.\\
34. Plagne, A. and Schmid, W.A., On the maximal cardinality of half-factorial sets in cyclic groups, Math. Ann. 333 (2005), no. 4, 759-785.\\
35. Radziejewski, M., A note on certain semigroups of algebraic numbers, Colloq. Math. 90 (2001), 51 - 58.\\
36. Radziejewski, M., On the distribution of algebraic numbers with prescribed factorization properties, Acta. Arith. 116 (2005), 153-171.\\
37. Radziejewski, M. and Schmid, W.A., On the asymptotic behavior of some counting functions, Colloq. Math. 102 (2005), 181-195.\\
38. Radziejewski, M., Schmid, W.A., Weakly half-factorial sets in finite abelian groups. Forum Math. 19 (2007), no. 4, 727-747.\\
39. Schmid, W.A., Arithmetic of block monoids, Math. Slovaca 54 (2004), 503-526.\\
40. Skula, L., On c-semigroups, Acta Arith. 31 (1976), 247-257.\\
41. Sliwa, J., On the length of factorizations in algebraic number fields, Bulletin de L'Academie Polonaise des Sciences, 23, No. 4(1975), { }387-389.
\\
42. Sliwa, J., A note on factorizations in algebraic number fields, Bull. Acad. Pol. Sci., Sr. Sci.Math. Astron. Phys. 24 (1976), 313 - 314.\\
43. Sliwa, J., Factorizations of distinct lengths in algebraic number fields, Acta Arith.31 (1976), 399-417.\\
44. Sliwa, J., Remarks on factorizations in algebraic number fields, Colloq. Math. 46 (1982), 123-130.\\
45. Zaks, A., Half factorial domains, Bull. Amer. Math. Soc. 82 (1976), 721-723.\\
\\

\end{document}